\documentclass[12pt]{iopart}

\usepackage{iopams}

\usepackage{cite}
\usepackage{graphicx}
\usepackage{bm}
\usepackage{color}

\begin{document}

\title{Numerical methods to evaluate Koopman matrix from system equations}

\author{Jun Ohkubo$^{1,2}$}
\address{$^1$ Graduate School of Science and Engineering, Saitama University, \\
255 Shimo-Okubo, Sakura-ku, Saitama 338-8570, Japan}
\address{$^2$ JST, PRESTO, 4-1-8 Honcho, Kawaguchi, Saitama 332-0012, Japan}
\ead{johkubo@mail.saitama-u.ac.jp}

\vspace{10pt}

\begin{abstract}
The Koopman operator is beneficial for analyzing nonlinear and stochastic dynamics; it is linear but infinite-dimensional, and it governs the evolution of observables. The extended dynamic mode decomposition (EDMD) is one of the famous methods in the Koopman operator approach. The EDMD employs a data set of snapshot pairs and a specific dictionary to evaluate an approximation for the Koopman operator, i.e., the Koopman matrix. In this study, we focus on stochastic differential equations, and a method to obtain the Koopman matrix is proposed. The proposed method does not need any data set, which employs the original system equations to evaluate some of the targeted elements of the Koopman matrix. The proposed method comprises combinatorics, an approximation of the resolvent, and extrapolations. Comparisons with the EDMD are performed for a noisy van der Pol system. The proposed method yields reasonable results even in cases wherein the EDMD exhibits a slow convergence behavior.
\end{abstract}

%
%
%
%
%

\section{Introduction}
\label{sec_introduction}

There are many works to predict and control dynamical and stochastic systems in the research fields in physics and engineering. Recently, the Koopman operator has attracted attention, which focuses on observables instead of the system states in dynamical and stochastic systems. Even when a system is nonlinear, the Koopman operator is linear albeit infinite-dimensional. Consequently, powerful tools used in operator theory, such as spectral analysis, can be employed to investigate nonlinear systems owing to the Koopman operator. After the significant works in \cite{Mezic2004,Mezic2005}, data-driven analysis methods such as dynamic mode decomposition \cite{Schmid2010} and extended dynamic mode decomposition (EDMD) \cite{Williams2015} were proposed and studied intensively. For example, with the aid of a data set of snapshot pairs and a dictionary, the EDMD yields the Koopman matrix, which corresponds to a finite-dimensional approximated representation of the Koopman operator. Then, the eigenvalues and eigenfunctions of the Koopman matrix possess rich information regarding the system dynamics; the information helps us understand the underlying properties of the systems. Hence, the data-driven approach of the Koopman operator formalism is used in various research fields such as spectral analysis \cite{Rowley2009,Kawahara2016,Takeishi2017,Li2017,Korda2020}, system control \cite{Korda2018,Li2020}, and model estimation \cite{Mauroy2016}. There are also applications to neural networks \cite{Dogra2020} and the analysis of algorithms \cite{Dietrich2020}. We can see overviews of the approach in \cite{Budisic2012, Mezic2013}, and the introductory part in the paper by \v{C}rnjari\'{c}-\v{Z}ic et al. \cite{Crnjaric-Zic2020} is a good guide for the various methods for the data-driven approach.

It is also possible to apply the Koopman operator approach to stochastic systems \cite{Williams2015, Klus2016, Crnjaric-Zic2020}. As discussed in \cite{Williams2015}, solutions of the backward Kolmogorov equation yield the Koopman operator. Spectral discussions for several stochastic differential equations were given in \cite{Crnjaric-Zic2020}; the data-driven approach can recover exact eigenfunctions well. There is another type of discussion for the Koopman operator; for example, the Koopman operator from the viewpoint of positive operator semigroups is denoted in \cite{Batkai_book}.

Although the Koopman operator can deal with statistics in stochastic differential equations, the duality in stochastic processes also yields them by solving different kinds of stochastic processes. The paper by Jansen and Kurt \cite{Jansen2014} is one of the good reviews for the duality in stochastic processes. For example, some types of stochastic differential equations have the corresponding birth-death processes; the solution of the dual birth-death process gives statistics for the original stochastic differential equation. We see examples in genetics in \cite{Carinci2015}. A systematic derivation for the dual process was discussed in \cite{Giardina2009}, and an extension of the duality was proposed in \cite{Ohkubo2013}. Although these works employed algebraic discussions, it was clarified that a simple derivation based on partial integration and basis expansions is enough to derive the dual-process \cite{Ohkubo2019}. Furthermore, the relations of orthogonal polynomials and duality were discussed in \cite{Ohkubo2019, Franceschini2019}. In addition to the mathematical discussions, there have been recent developments for numerical methods for the applications of duality. Although a naive numerical evaluation based on the Monte Carlo methods takes a long computational time, a numerical method based on combinatorics and resolvent was proposed in \cite{Ohkubo2021}, which enables us to obtain statistics of stochastic differential equations in reasonable computational times without using any sampling method.

Here, let us consider the following question: Is it possible to obtain the Koopman matrix directly from system equations? The answer is yes; it is enough to discuss numerical solutions of the backward Kolmogorov equations or their generators. However, sometimes such approaches could become computationally challenging; if the dimension of the system is high, the computational cost becomes large. Then, it would be helpful if the numerical method developed in the duality is available to evaluate the Koopman matrix.

In the present paper, a numerical method to evaluate the Koopman matrix for stochastic differential equations is proposed. The proposed method employs the idea based on combinatorics and resolvent developed in \cite{Ohkubo2021}. The key point is as follows: we focus on only some elements of the Koopman matrix instead of the whole one. Hence, we expect a reduction in computational costs. The connection of the duality with the Koopman operator is clarified, and an extrapolation technique is employed to reduce the computational time. Note that the relationship between the duality and the Koopman approach in control theory was discussed in \cite{Ohkubo2021a}. However, the discussion was a naive one, and it did not intend the computation of the Koopman matrix. The present paper gives the first attempt to compute the Koopman matrix elements. Although much mathematically rigorous discussions remain to be done, the numerical demonstrations here will become motivations for future mathematical and practical works.

The remainder of the paper is organized as follows. In section~2, the Koopman operator and the EDMD method are reviewed. Section~3 presents the primary contribution of this study, i.e., the connection of the duality with the Koopman operator formalism and an explicit algorithm based on the resolvent and extrapolation. Section~4 gives numerical results of the proposed method; a noisy van der Pol system is used as a concrete example. Moreover, comparisons with the EDMD method are also given. Finally, section~5 concludes this paper.

\section{Previous knowledge for Koopman operator}

\subsection{Koopman operator}
\label{subsec_Koopman}

First, a brief explanation for the Koopman operator is provided. For further details, see \cite{Mauroy2020}.

Consider $\mathcal{M} \subseteq \mathbb{R}^{D}$ to be a state space. For a state $\bm{x} \in \mathcal{M}$, consider a discrete-time deterministic process
\begin{eqnarray}
\bm{x} \mapsto \bm{F}(\bm{x}),
\end{eqnarray}
where $\bm{F}: \mathcal{M} \to \mathcal{M}$ is the time-evolution operator. That is, for a discrete time step $\ell \in \mathbb{Z}$, a dynamical system is defined as follows:
\begin{eqnarray}
\bm{x}(\ell+1) = \bm{F}(\bm{x}(\ell)), \quad \bm{x}(\ell) \in \mathcal{M}.
\end{eqnarray}
The Koopman formalism focuses on the time-evolution of observables, which are represented by functions on the state space $\mathcal{M}$ in a suitable function space as follows. Consider the Hilbert space
\begin{eqnarray}
L^2(\mathcal{M},\rho) = \left\{ \phi: \mathcal{M} \to \mathbb{C} : \| \phi \|_{L^2(\mathcal{M},\rho)} < \infty \right\},
\end{eqnarray}
where the inner product is defined as
\begin{eqnarray}
\left( \phi_1, \phi_2 \right) = \int_\mathcal{M} \phi_1(\bm{x}) \overline{\phi_2(\bm{x})} \rho(\rmd\bm{x})
\end{eqnarray}
and the norm is expressed as
\begin{eqnarray}
\| \phi \|_{L^2(\mathcal{M},\rho)}
= \left( 
\int_\mathcal{M} \left| \phi(\bm{x}) \right|^2 \rho(\rmd\bm{x})
\right)^{\frac{1}{2}}.
\end{eqnarray}
Consequently, the Koopman operator $\mathcal{K}$ is defined as an operator which acts on an observable $\phi \in L^2(\mathcal{M},\rho)$ as follows:
\begin{eqnarray}
\mathcal{K} \phi = \phi \circ \bm{F}.
\end{eqnarray}
The above settings for the observables are the same as \cite{Williams2015}. Although the $L^2$ condition could be relaxed, the discussion in \cite{Williams2015} needed the inner products in the Galerkin-like method. Hence, the $L^2$ condition was employed. In the present paper, to connect the Koopman approach with the duality in the stochastic process, we use the same $L^2$ condition. In addition, it is suitable to assume that observables are expressed as a linear combination of orthogonal polynomials such as Hermite polynomials. Then, $\rho$ is an invariant measure, and it is naturally determined from the orthogonal polynomials. These settings are enough for the following discussions. (Of course, the more mathematically rigorous discussion may clarify that a more relaxed condition is possible; this is beyond the scope of the present paper.)

As mentioned in section~1, there is a stochastic version for the Koopman operator. For a discrete-time Markov process with
\begin{eqnarray}
\bm{x} \mapsto \bm{F}(\bm{x}; \bm{\omega}),
\end{eqnarray}
the stochastic Koopman operator is defined as
\begin{eqnarray}
(\mathcal{K} \phi)(\bm{x}) = \mathbb{E}\left[ \phi \left( \bm{F}(\bm{x}; \bm{\omega}) \right)\right],
\end{eqnarray}
where $\bm{\omega} \in \Omega_\mathrm{s}$ is an element in the probability space associated with the stochastic dynamics $\Omega_\mathrm{s}$ and the probability measure $P$. Further, $\mathbb{E}[\cdot]$ denotes the expected value over that space, and $\phi \in L^2(\mathcal{M},\rho)$ is a scalar observable. The references \cite{Williams2015, Crnjaric-Zic2020} include enough information for the stochastic cases.

The Koopman operator is linear even when the system exhibits nonlinearity. Although the Koopman operator is infinite-dimensional, the linear action of the Koopman operator is preferable from the perspective of computation. In practice, the Koopman operator is approximated by the corresponding finite-dimensional matrix representation. This characteristic of linearity is useful in spectral analysis and mode analysis.

The remaining task here is to obtain the Koopman operator numerically. The Koopman matrix is usually employed to approximate the Koopman operator. In addition, there are several works for data-driven methods to evaluate the Koopman matrix numerically. The following section reviews one of the methods.

\subsection{Extended Dynamic Mode Decomposition}
\label{subsec_EDMD}

A finite-dimensional matrix cannot represent the Koopman operator exactly. Thus, a certain approximation is necessary. Consequently, the Koopman matrix is obtained instead of the Koopman operator. The Koopman matrix has a finite dimension, and its dimension depends on the size of the dictionary. Extended dynamic mode decomposition (EDMD) is one of the methods to evaluate the Koopman matrix from data \cite{Williams2015}.

The EDMD requires a data set of snapshot pairs and a dictionary of observables. Let $N_\mathrm{data}$ be the number of snapshot pairs. The snapshot pairs are denoted as $\{ (\bm{x}_n, \bm{y}_n) \}_{n=1}^{N_\mathrm{data}}$, where 
\begin{eqnarray}
\bm{x}_n = 
(
\begin{array}{cccc}
x_{n1} \quad x_{n2} \quad \cdots \quad x_{nD}
\end{array}
)^\mathrm{T}
\end{eqnarray}
and 
\begin{eqnarray}
\bm{y}_n = \bm{F}(\bm{x}_n).
\label{eq_discrete_time_evolution}
\end{eqnarray}
Approximating the Koopman operator $\mathcal{K}$ in the form of a finite-dimensional matrix $K$ requires the dictionary of observables. The dictionary $\bm{\psi}(\bm{x})$ with the size $N_\mathrm{dic}$ is denoted as
\begin{eqnarray}
\bm{\psi}(\bm{x}) 
= \left(
\begin{array}{cccc}
\psi_1(\bm{x}) \quad \psi_2(\bm{x}) \quad \cdots \quad \psi_{N_\mathrm{dic}}(\bm{x})
\end{array}
\right)^\mathrm{T},
\end{eqnarray}
where $\psi_i: \mathcal{M} \to \mathbb{C}$ for $i = 1,\dots,N_\mathrm{dic}$. The dictionary is chosen to approximate the Koopman operator well, which spans a subspace $\mathcal{F}_{\bm{\psi}} \subseteq L^2(\mathcal{M},\rho)$. Using the dictionary, an observable $\phi(\bm{x}) \in L^2(\mathcal{M},\rho)$ can be represented as
\begin{eqnarray}
\left\{ 
\begin{array}{ll}
\phi(\bm{x}) = \sum_{i=1}^{N_\mathrm{dic}} a^{\phi}_i \psi_i (\bm{x}) & \textrm{for} \,\, \phi(\bm{x}) \in \mathcal{F}_{\bm{\psi}},\\
\phi(\bm{x}) \simeq \sum_{i=1}^{N_\mathrm{dic}} a^{\phi}_i \psi_i (\bm{x}) & \textrm{for} \,\, \phi(\bm{x}) \notin \mathcal{F}_{\bm{\psi}},
\end{array}
\right.
\label{eq_function_basis_expansion}
\end{eqnarray}
where $\{a^{\phi}_i\}_{i=1}^{N_\mathrm{dic}}$ are the expansion coefficients. The finite-size dictionary is employed here, and hence observables outside of the subspace $\mathcal{F}_{\bm{\psi}}$ are only approximately represented. Then, the action of the Koopman operator is expressed as
\begin{eqnarray}
\mathcal{K} \phi(\bm{x}) \simeq \sum_{i=1}^{N_\mathrm{dic}} a^{\phi}_i \left( \mathcal{K} \psi_i (\bm{x}) \right).
\label{eq_action_Koopman}
\end{eqnarray}
If $\mathcal{K}\phi(\bm{x}) \in \mathcal{F}_{\bm{\psi}}$, \eref{eq_action_Koopman} yields an exact equality.

With the aid of the dictionary $\bm{\psi}(\bm{x})$, the Koopman operator is approximated by a finite-size matrix $K$ as follows:
\begin{eqnarray}
\left(\begin{array}{c}
\mathcal{K} \psi_1(\bm{x}) \\ \mathcal{K} \psi_2(\bm{x}) \\ \vdots \\ \mathcal{K} \psi_{N_\mathrm{dic}}(\bm{x}) 
\end{array}\right)
=
\mathcal{K} \bm{\psi}(\bm{x})
\simeq
K^{\mathrm{T}} \bm{\psi}(\bm{x}).
\label{eq_action_Koopman_operator_and_matrix}
\end{eqnarray}
To determine the Koopman matrix $K$, the following cost function is introduced:
\begin{eqnarray}
J(K) = \sum_{n=1}^{N_\mathrm{data}} \left\| \bm{\psi}(\bm{y}_n) - K^{\mathrm{T}} \bm{\psi}(\bm{x}_n) \right\|^2,
\label{eq_Koopman_cost_function}
\end{eqnarray}
which focuses on the consistency of the time-evolution of the dictionary for the snapshot pairs. 
By minimizing the cost function \eref{eq_Koopman_cost_function}, we obtain the following result:
\begin{eqnarray}
K^{\mathrm{T}} = A G^{+},
\label{eq_EDMD_algorithm}
\end{eqnarray}
where 
\begin{eqnarray}
A &= \frac{1}{N_\mathrm{data}} \sum_{n=1}^{N_\mathrm{data}} \bm{\psi}(\bm{y}_n) \bm{\psi}(\bm{x}_n)^\mathrm{T}, \\
G &= \frac{1}{N_\mathrm{data}} \sum_{n=1}^{N_\mathrm{data}} \bm{\psi}(\bm{x}_n) \bm{\psi}(\bm{x}_n)^\mathrm{T},
\end{eqnarray}
and $G^{+}$ is the pseudo-inverse of $G$. As for the details of the derivation, for example, see \cite{Klus2016}.

Once we obtain the Koopman matrix $K$, it is possible to evaluate the values of observables following the time-evolution with the linear action of the matrix $K^\mathrm{T}$. That is, 
\begin{eqnarray}
\phi(\bm{y}_n) \simeq 
(\begin{array}{cccc}
a_1^{\phi} \quad a_2^{\phi} \quad \cdots \quad a_{N_\mathrm{dic}}^{\phi}
\end{array})
\,\,  K^{\mathrm{T}}
\left(\begin{array}{c}
\psi_1(\bm{x}_n) \\ \psi_2(\bm{x}_n) \\ \vdots \\ \psi_{N_\mathrm{dic}}(\bm{x}_n)
\end{array}\right).
\end{eqnarray}

Usually, it is possible to choose an arbitrary dictionary. However, the following discussion employs a monomial-type dictionary. Let $\bm{\alpha}_s = (\alpha_{1,s}, \dots, \alpha_{D,s})$ and $\alpha_{d,s} \in \mathbb{N}_0$ for $d = 1, \dots, D$. Then, 
\begin{eqnarray}
\psi_{s}(\bm{x}) = x_1^{\alpha_{1,s}} \cdots x_D^{\alpha_{D,s}} = \bm{x}^{\bm{\alpha}_s}.
\label{eq_monomial_dictionary}
\end{eqnarray}
This choice of the dictionary naturally connects the Koopman matrix evaluated from the EDMD to the duality relation for stochastic differential equations. Hence, in the subsequent sections, the dictionary with \eref{eq_monomial_dictionary} is employed.

There is a comment for the definition of the Koopman matrix. The usage of the transpose matrix as in \eref{eq_action_Koopman_operator_and_matrix} is conventional in the EDMD. As pointed out in \cite{Klus2016}, this definition means that the matrix representation of Ulam's method has the same form as the Koopman matrix. Since Ulam's method is the dual one of the Koopman approach, one may consider that the Koopman matrix should be defined as its transposed one. In the following discussion, the Koopman matrix is connected to the duality in the stochastic process, and the matrix element directly corresponds to transitions from an initial state to a final state. Hence, the following definition is useful:
\begin{eqnarray}
\widetilde{K} \equiv K^\mathrm{T}.
\end{eqnarray}
In the following discussions, we employ this definition $\widetilde{K}$ as the Koopman matrix.

\subsection{Koopman operator for stochastic differential equations}
\label{subsec_sde}

In this work, we focus only on stochastic differential equations. Here, stochastic differential equations and their Koopman operator are reviewed concisely. For details, see \cite{Crnjaric-Zic2020}. 

Consider a state space $\mathcal{M} \subseteq \mathbb{R}^D$ and a vector of stochastic variables $\bm{X} \in \mathcal{M}$. In contrast to section~\ref{subsec_Koopman}, herein, a system is treated with a continuous time $t\in \mathbb{R}$. The time-evolution of the system is provided as the following coupled stochastic differential equations \cite{Gardiner_book}:
\begin{eqnarray}
\rmd\bm{X} = \bm{a}(\bm{X}, t) \rmd t + B(\bm{X}, t) \rmd\bm{W}(t),
\label{eq_SDE}
\end{eqnarray}
where $\bm{a}(\bm{X}, t)$ is a vector of drift coefficient functions and $B(\bm{X}, t)$ is a matrix of diffusion coefficient functions. Further, components of the vector of Wiener processes $\bm{W}(t)$ satisfy
\begin{eqnarray}
\rmd W_i(t) \rmd W_j(t) = \delta_{ij} \rmd t.
\end{eqnarray}
The stochastic differential equations in \eref{eq_SDE} have a corresponding Fokker-Planck equation as follows \cite{Gardiner_book}:
\begin{eqnarray}
\fl
\frac{\partial}{\partial t} p(\bm{x},t)
= - \sum_{i} \frac{\partial}{\partial x_i} \left( a_i(\bm{x},t) p(\bm{x},t) \right)
+ \frac{1}{2} \sum_{i,j} \frac{\partial^2}{\partial x_i \partial x_j} 
\left( \left[ B(\bm{x},t) B(\bm{x},t)^\mathrm{T} \right]_{ij} p(\bm{x},t) \right),
\label{eq_Fokker-Planck}
\end{eqnarray}
where $p(\bm{x},t)$ is the probability density function for $\bm{x}\in\mathcal{M}$ at time $t$. When we consider the time-evolution from $t_\mathrm{i}$ to $t_\mathrm{f}$, an expectation value for $\phi \in L^2(\mathcal{M},\rho)$ is evaluated as
\begin{eqnarray}
\mathbb{E}_{t_\mathrm{f}}\left[ \phi(\bm{X}) \right] 
= \int_{\mathcal{M}} \phi(\bm{x}) p(\bm{x},t_\mathrm{f}) \rmd\bm{x}.
\label{eq_expectation_basic}
\end{eqnarray}

In general, it is a challenging task to obtain an explicit solution for the probability density. Hence, Monte Carlo simulations for the stochastic differential equations are usually employed to generate sample paths, and an approximate solution for $p(\bm{x}, t_\mathrm{f})$ is numerically estimated. However, the Koopman formalism enables us to obtain the expectation values of statistics, \eref{eq_expectation_basic}, without solving the Fokker-Planck equation and using the Monte Carlo simulations; the stochastic Koopman operator $\mathcal{K}$ yields the corresponding function $\varphi \equiv \mathcal{K} \phi$ from the observable $\phi$. Here, it is assumed that an initial condition for the probability density function at time $t_\mathrm{i}$ is
\begin{eqnarray}
p(\bm{x}, t_\mathrm{i}) = \delta(\bm{x} - \bm{x}_\mathrm{i}),
\end{eqnarray}
where $\delta(\cdot)$ is the Dirac delta function, and $\bm{x}_\mathrm{i}$ is the initial position. By focusing only on the final results of the time-evolution of the stochastic differential equations (or the corresponding Fokker-Planck equation) from $t_\mathrm{i}$ to $t_\mathrm{f}$, it is possible to consider the time-evolution as a discrete-time system. Consequently, there exists a corresponding Koopman operator $\mathcal{K}$:
\begin{eqnarray}
\left( \mathcal{K}\phi\right) (\bm{x}_{\mathrm{i}}) 
= \varphi(\bm{x}_{\mathrm{i}}) 
= \mathbb{E}_{t_\mathrm{f}} \left[ \phi(\bm{X})\right].
\end{eqnarray}
It is possible to employ the EDMD in section~\ref{subsec_EDMD} to obtain the corresponding Koopman matrix for the stochastic differential equations.

\section{Proposal of numerical methods for evaluating Koopman matrix}
\label{sec_proposal}

This section presents the numerical method to obtain the Koopman matrix for stochastic differential equations, which is the primary contribution of the present paper. Based on discussions for the duality in the stochastic process, the proposed method utilizes the system equations, combinatorics, and resolvents of the time-evolution operator. In the proposed method, an approximation for the resolvent is employed. Furthermore, an extrapolation technique is combined for efficient numerical evaluation. Then, the proposed algorithm directly computes a matrix element of the Koopman matrix.

\subsection{Requirements of the stochastic differential equations}

First of all, a requirement for the proposed method should be mentioned. As we will see soon, the following discussion needs the condition that drift and diffusion coefficient functions should be represented as polynomials. Although one may consider that the requirement is a little tight, an extension of stochastic variables enables us to apply the proposed method for various types of systems. This extension technique was proposed in [1] in the context of the duality in the stochastic process. 

For example, if there is a term of $\sin(X)$ in a coefficient function, it is enough to define the following two variables:
\begin{eqnarray}
Z_1 = \sin(X), \\
Z_2 = \cos(X).
\end{eqnarray}
Since $X$ is a stochastic variable, $Z_1$ and $Z_2$ are also stochastic variables. Note that such additional variables are expressed as functions for $X$. Hence, The following Ito formula gives the corresponding stochastic differential equations for the newly introduced variables \cite{Gardiner_book}:
\begin{eqnarray}
\fl
\rmd f(\bm{X}) =& 
\left\{
\sum_i a_i(\bm{X},t) \partial_{X_i} f(\bm{X})
+ \frac{1}{2} \sum_{i,j} \left[ B(\bm{X},t) B^\mathrm{T}(\bm{X},t) \right]_{ij} \partial_{X_i} \partial_{X_j} f(\bm{X})
\right\} \rmd t \nonumber \\
\fl
&+ \sum_{i,j} B_{ij}(\bm{X},t) \partial_{X_i} f(\bm{X}) \, \rmd W_j(t),
\end{eqnarray}
where $f(\bm{X})$ is a function that we wish to replace. Furthermore, the replacement of the time variable $t$ in such functions as a new variable removes the explicit time dependency in $\{a_i(\bm{X},t)\}$ and $B(\bm{X},t)$. Employing the above technique, functions such as $\exp$ and $\log$ are also available. The mathematically rigorous condition for what types of functions we can use is still unclear. This feature may be related to the concept of holonomic functions \cite{Takayama2013}, and we need further mathematical studies in the future. However, we can expect that most practical problems could satisfy this feature. 

As a consequence of the above discussions, we assume the polynomial types of coefficient functions in the following discussions. Note that the dimension of the extended system is different from the original one for the state space in section~\ref{subsec_sde}. However, there is no problem when we evaluate explicit values of the Koopman matrix elements; the newly introduced state variables are evaluated numerically with concrete initial values, as shown below. Hence, we write the extended variables simply as $\bm{X}$ or $\bm{x}$.

\subsection{Evaluation of Koopman matrix from dual process}
\label{subsec_dual}

Firstly, the Fokker-Planck equation is rewritten in \eref{eq_Fokker-Planck} as follows:
\begin{eqnarray}
\frac{\partial}{\partial t} p(\bm{x},t) = \mathcal{L} \, p(\bm{x},t),
\end{eqnarray}
where
\begin{eqnarray}
\mathcal{L}
= - \sum_{i} \frac{\partial}{\partial x_i} a_i(\bm{x},t) 
+ \frac{1}{2} \sum_{i,j} \frac{\partial^2}{\partial x_i \partial x_j} 
\left[ B(\bm{x},t) B(\bm{x},t)^\mathrm{T} \right]_{ij}
\end{eqnarray}
is the time-evolution operator. Assume that the initial condition of the Fokker-Planck equation obeys the Dirac delta function. Then, it is possible to calculate expectation values as follows:
\begin{eqnarray}
\mathbb{E}\left[ x_1^{\alpha_1} \cdots x_D^{\alpha_D} \right]
& = \mathbb{E}\left[ \bm{x}^{\bm{\alpha}} \right] \nonumber \\
&= \int \bm{x}^{\bm{\alpha}} p(\bm{x},t) \rmd\bm{x} \nonumber \\
&= \int \bm{x}^{\bm{\alpha}} \left( e^{\mathcal{L} t} \delta(\bm{x}-\bm{x}_0)\right) \rmd\bm{x} \nonumber \\
&= \int \left( e^{\mathcal{L}^\dagger t} \bm{x}^{\bm{\alpha}} \right) \delta(\bm{x}-\bm{x}_0) \rmd\bm{x} \nonumber \\
&= \int \varphi_{\bm{\alpha}}(\bm{x},t) \delta(\bm{x}-\bm{x}_0) \rmd\bm{x} \nonumber \\
&= \varphi_{\bm{\alpha}}(\bm{x}_0,t),
\label{eq_derivation_expectation}
\end{eqnarray}
where
\begin{eqnarray}
\mathcal{L}^\dagger
= \sum_{i} a_i(\bm{x},t) \frac{\partial}{\partial x_i} 
+ \frac{1}{2} \sum_{i,j} \left[ B(\bm{x},t) B(\bm{x},t)^\mathrm{T} \right]_{ij} \frac{\partial^2}{\partial x_i \partial x_j} 
\label{eq_adjoint_operator}
\end{eqnarray}
is the adjoint operator of $\mathcal{L}$. As discussed in \cite{Ohkubo2019}, it is easy to derive the adjoint operator by performing integration by parts. Further, function $\varphi_{\bm{\alpha}}(\bm{x},t)$ obeys the following time-evolution equation (the backward Kolmogorov equation):
\begin{eqnarray}
\frac{d}{dt} \varphi_{\bm{\alpha}}(\bm{x},t) = \mathcal{L}^\dagger \varphi_{\bm{\alpha}}(\bm{x},t),
\label{eq_time_evolution_adjoint}
\end{eqnarray}
where the initial condition is $\varphi_{\bm{\alpha}}(\bm{x},t_\mathrm{i}) = \bm{x}^{\bm{\alpha}}$. Note that $\varphi_{\bm{\alpha}}(\bm{x},t)$ is not a probability density function; there is no probability conservation law.

Next, consider the following expansion for $\varphi_{\bm{\alpha}}(\bm{x},t)$:
\begin{eqnarray}
\varphi_{\bm{\alpha}}(\bm{x},t) = \sum_{\bm{n}} P_{\bm{\alpha}}(\bm{n},t) \bm{x}^{\bm{n}},
\label{eq_basis_expansion_1}
\end{eqnarray}
where $\bm{n} = (n_1, \dots, n_D)$ and $n_d \in \mathbb{N}_0$. $\{P_{\bm{\alpha}}(\bm{n},t)\}$ are the expansion coefficients. Here, the simple functions $\{\bm{x}^{\bm{n}}\}$ are employed as the basis.

There is a comment for the choice of the basis. The basis corresponds to the dictionary in the EDMD, as shown below. Of course, it is possible to use other basis functions; in the context of the duality, Hermite polynomials were used in \cite{Ohkubo2019}, and Legendre polynomials were employed to deal with variables with finite ranges in \cite{Ohkubo2020}. For simplicity, the following explanation is restricted to the cases with the simple basis in \eref{eq_basis_expansion_1}.

Using the ket-notation $\left| \bm{n} \right\rangle \equiv \bm{x}^{\bm{n}}$, \eref{eq_basis_expansion_1} can be rewritten as
\begin{eqnarray}
\varphi_{\bm{\alpha}}(\bm{x},t) = \sum_{\bm{n}} P_{\bm{\alpha}}(\bm{n},t) \left| \bm{n} \right\rangle.
\end{eqnarray}
When the coefficient functions in \eref{eq_adjoint_operator} are expressed only with simple polynomials of $\bm{x}$, the action of $\mathcal{L}^\dagger$ on $\varphi_{\bm{\alpha}}(\bm{x},t)$ yields coupled ordinary differential equations for $\{P_{\bm{\alpha}}(\bm{n},t)\}$. As shown later, the coupled ordinary differential equations give the dual process. It is easy to see that the coefficients $\{P_{\bm{\alpha}}(\bm{n},t)\}$ directly correspond to the elements of the Koopman matrix, with $\bm{\alpha}$ and $\bm{n}$ indicating the row and column of the matrix, respectively.

As reported in \cite{Ohkubo2013,Ohkubo2019}, it is possible to recover stochasticity for \eref{eq_time_evolution_adjoint}. Then, the derived equations correspond to chemical master equations, which can be numerically solved using conventional Runge-Kutta methods or Monte Carlo methods such as the Gillespie algorithm \cite{Gillespie1977}. However, the approach based on coupled ordinary differential equations is not employed here. Sometimes the direct time-evolution of the coupled ordinary differential equations is intractable because of the high dimensionality. Moreover, sometimes there could be no need to evaluate all the elements of the Koopman matrix at once. The subsequent discussion introduces the method based on combinatorics. The proposed method is suitable to compute a targeted matrix element.

As assumed above, coefficient functions in the time-evolution operator $\mathcal{L}^\dagger$ in \eref{eq_adjoint_operator} have simple polynomial forms. That is, both $a_i(\bm{x},t)$ and $B(\bm{x},t)$ are expressed as polynomials in $\bm{x}$. Consequently, the action of each term in $\mathcal{L}^\dagger$ simply changes the state vector $|\bm{n}\rangle$. For example, 
\begin{eqnarray}
\fl
x_d | \bm{n} \rangle
= x_d \times \left( x_1^{n_1} \dots x_d^{n_d} \dots x_D^{n_D} \right)
= x_1^{n_1} \dots x_d^{n_d+1} \dots x_D^{n_D} 
= \left| \bm{n} + \bm{1}_d \right\rangle,
\end{eqnarray}
where $\bm{1}_d$ is a standard basis vector whose components are all zero, except $d$-th one that equals $1$. In a similar manner,
\begin{eqnarray}
\fl
\partial_{x_d} | \bm{n} \rangle
= \partial_{x_d} \left( x_1^{n_1} \dots x_d^{n_d} \dots x_D^{n_D} \right)
= n_d \, x_1^{n_1} \dots x_d^{n_d-1} \dots x_D^{n_D} 
= n_d \left| \bm{n} - \bm{1}_d \right\rangle.
\end{eqnarray}
In general, the state vector is represented as the direct product
\begin{eqnarray}
|\bm{n}\rangle = | n_1 \rangle \otimes \dots \otimes | n_D \rangle,
\end{eqnarray} 
and it is possible to express the action of a term in $\mathcal{L}^\dagger$ by using the following form:
\begin{eqnarray}
x_d^{m^+} \frac{\partial^{m^-}}{\partial x_d^{m^-}} \left| n_d \right\rangle 
&= x_d^{m^+} \frac{\partial^{m^-}}{\partial x_d^{m^-}} x_d^{n_d} \nonumber \\
&= \frac{n_d!}{(n_d-m^-)!} x_d^{n_d-m^-+m^+} \nonumber \\
&= \frac{n_d!}{(n_d-m^-)!} \left| n_d-m^-+m^+ \right\rangle
\label{eq_general_action}
\end{eqnarray}
for $m^- \le n_d$; otherwise, the r.h.s in \eref{eq_general_action} is zero.

Here, the `bra' notation is defined as follows:
\begin{eqnarray}
\langle \bm{n}' | = \int \rmd \bm{x} \, \delta(\bm{x}) \left(\frac{\rmd}{\rmd \bm{x}} \right)^{\bm{n}'} (\cdot).
\end{eqnarray}
The bra and ket vectors satisfy 
\begin{eqnarray}
\langle \bm{n}' | \bm{n} \rangle = \bm{n}! \, \delta_{\bm{n}',\bm{n}},
\label{eq_bra_ket_orthogonal}
\end{eqnarray}
where $\bm{n}! \equiv n_1! n_2! \cdots n_D!$. Then, we obtain the following decomposition of unity:
\begin{eqnarray}
\bm{1} = \sum_{\bm{n}} \frac{1}{\bm{n}!} | \bm{n} \rangle \langle \bm{n} |.
\label{eq_unity}
\end{eqnarray}
It might be difficult to image the function space defined by the above bra and ket notations. As stated above, the Hermite polynomials can give essentially the same formalism, and linear combinations of the Hermite polynomials yield the monomial-type basis. Hence, it might be better to restrict the discussion within the Hermite cases. However, the monomial-type basis is suitable for the following derivation. Furthermore, numerical results for some examples show that empirically the simple basis works well. Hence, we employ these definitions and avoid further mathematical rigorous discussions.

Here, it is convenient to treat the drift and diffusion terms in \eref{eq_adjoint_operator} without distinguishing between them. After expanding $\mathcal{L}^\dagger$, let $R$ be the number of terms in $\mathcal{L}^\dagger$ in the expanded form. Then, acting $\mathcal{L}^\dagger$ to the decomposition of unity in \eref{eq_unity} from the left, $\mathcal{L}^\dagger$ in \eref{eq_adjoint_operator} is rewritten as the following notation:
\begin{eqnarray}
\mathcal{L}^\dagger = \sum_{r=1}^R \sum_{\bm{n}} \gamma_r(\bm{n}) \frac{1}{\bm{n}!}\, | \bm{n} + \bm{v}_r \rangle \langle \bm{n} |,
\label{eq_simple_time_evolution_operator}
\end{eqnarray}
where $\gamma_r(\bm{n})$ is the state-dependent coefficient, and $\bm{v}_r$ is the vector for the state change, which can be calculated using \eref{eq_general_action}. In this study, the target is the stochastic differential equations. Hence, functional forms of the terms in $\mathcal{L}^\dagger$ have been restricted into two cases. The first case is as follows:
\begin{eqnarray}
c \, x_1^{m_1} \cdots x_D^{m_D} \frac{\partial}{\partial x_i} | \bm{n} \rangle
= c \, n_i | \bm{n} + \bm{v}_i \rangle,
\label{eq_gamma_1}
\end{eqnarray}
where $c$ is a state-independent constant, $m_i \in \mathbb{N}_0$ for $i = 1, \dots, D$, and
\begin{eqnarray}
\bm{v}_i 
= \left[ m_1, \dots, m_i-1, \dots, m_D \right].
\end{eqnarray}
The second case is
\begin{eqnarray}
c \, x_1^{m_1} \cdots x_D^{m_D} \frac{\partial}{\partial x_i} \frac{\partial}{\partial x_j} \left| \bm{n} \right\rangle
= \left\{ \begin{array}{ll}
c \, n_i (n_i-1) \left| \bm{n} + \bm{v}_i \right\rangle & \textrm{if $i = j$,}\\
c \, n_i n_j \left| \bm{n} + \bm{v}_i \right\rangle & \textrm{otherwise,}
\end{array}\right.
\label{eq_gamma_2}
\end{eqnarray}
where
\begin{eqnarray}
\fl
\bm{v}_i 
= \left\{ \begin{array}{ll}
\left[ m_1, \dots, m_i-2, \dots, m_D \right] & \textrm{if $i = j$,}\\
 & \\
\left[ m_1, \dots, m_i-1, \dots, m_j-1, \dots, m_D \right] & \textrm{otherwise.}
\end{array}\right.
\end{eqnarray}
Using \eref{eq_gamma_1} and \eref{eq_gamma_2}, the coefficient functions $\{\gamma_r(\bm{n})\}$ in \eref{eq_simple_time_evolution_operator} can be explicitly calculated.

Let $T$ be the time-interval of the time-evolution, that is, $T \equiv t_\mathrm{f} - t_\mathrm{i}$. Consequently, the solution is formally written as:
\begin{eqnarray}
\varphi_{\alpha}(\bm{x},T) = e^{\mathcal{L}^\dagger T} | \bm{\alpha} \rangle.
\end{eqnarray}
To evaluate $P_{\bm{\alpha}}(\bm{n},T)$, the time-evolution operator is split with a small-time interval $\Delta t$, that is, 
\begin{eqnarray}
e^{\mathcal{L}^\dagger T} = e^{\mathcal{L}^\dagger \Delta t} \dots e^{\mathcal{L}^\dagger \Delta t}.
\label{eq_time_evolution_split}
\end{eqnarray}
Subsequently, the Taylor expansion provides 
\begin{eqnarray}
e^{\mathcal{L}^\dagger \Delta t} 
&\simeq 
1 + \sum_{r=1}^R \sum_{\bm{n}} \gamma_r(\bm{n}) \frac{1}{\bm{n}!} \Delta t | \bm{n} + \bm{v}_r \rangle \langle \bm{n} |.
\label{eq_adjoint_expansion_reinterpretation}
\end{eqnarray}
As discussed earlier, the action of the time-evolution operator results in a change in the discrete state vector. Hence, the \textit{dual} process of the original stochastic differential equation is derived. As stated above, it is possible to recover the stochastic property, which results in the corresponding dual stochastic process. However, there is no need to recover stochasticity in the following discussions.

The repeated insertion of the decomposition of unity in \eref{eq_unity} into the time-evolution operator in \eref{eq_time_evolution_split} leads to the following path-integral-like result for the initial condition $|\bm{n}_\mathrm{ini} \rangle$:
\begin{eqnarray}
\fl
e^{\mathcal{L}^\dagger T} | \bm{n}_{\mathrm{ini}} \rangle 
\simeq
\sum_{M=0}^\infty \frac{T^M}{M!} \left( \sum_{j^{(1)}} \dots \sum_{j^{(M)}}
\gamma_{j^{(M)}}(\bm{n}_{M-1}) \cdots \gamma_{j^{(1)}}(\bm{n}_{\mathrm{ini}}) 
 \left| \bm{n}_{\mathrm{ini}} + \sum_{m=1}^M \bm{v}_{j^{(m)}} \right\rangle \right), \nonumber \\
\fl
\label{eq_path_integral_like}
\end{eqnarray}
where the expansion in \eref{eq_adjoint_expansion_reinterpretation} is employed; a supplementary explanation for \eref{eq_path_integral_like} is given in \ref{sec_appendix_path_integral}. Note that $(1/\bm{n}!) \langle \bm{n}| \exp(\mathcal{L}^\dagger T) | \bm{n}_\mathrm{ini} \rangle$ directly corresponds to $P_{\bm{n}_\mathrm{ini}} (\bm{n},t)$. As discussed above, $P_{\bm{n}_\mathrm{ini}} (\bm{n},t)$ is equal to the corresponding element of the Koopman matrix. However, obtaining the summation for all the paths in \eref{eq_path_integral_like} is challenging; that is, the $M$ divisions should be sufficiently large, and hence the computational complexity is high. In addition, as shown in \cite{Ohkubo2021}, the summation in \eref{eq_path_integral_like} exhibits a diverging behavior. Therefore, a different calculation method should be employed to evaluate the exponential of the time-evolution operator.

\subsection{Resolvent and extrapolation technique}
\label{subsec_resolvent}

The exponential of the time-evolution operator is naively expanded as follows:
\begin{eqnarray}
\rme^{\mathcal{L}^\dagger T}
= \sum_{m=0}^\infty \frac{T^m}{m!} (\mathcal{L}^\dagger)^m.
\end{eqnarray}
The definition based on the Taylor series is numerically unusable if $\mathcal{L}^\dagger$ is an unbounded operator in a Banach space \cite{Kato_book}. Alternatively, the following formula is suitable for the numerical calculation:
\begin{eqnarray}
\rme^{\mathcal{L}^\dagger T}
= \lim_{M \to \infty} \left[ \left(1 - \frac{T}{M} \mathcal{L}^\dagger \right)^{-1} \right]^M,
\label{eq_resolvent}
\end{eqnarray}
where $\left(1 - \frac{T}{M} \mathcal{L}^\dagger \right)^{-1}$ is a resolvent of $\mathcal{L}^\dagger$, apart from a constant factor \cite{Kato_book}.

Here, the conventional Gauss–Jordan elimination is employed to evaluate the resolvent approximately \cite{Strang_book}. It is possible to rewrite the element of the (infinite-dimensional) matrix $\mathcal{L}^\dagger$ as
\begin{eqnarray}
\fl
\left[ \mathcal{L}^\dagger \right]_{s's}
= \frac{1}{\bm{n}_{s'}!} \langle \bm{n}_{s'} | 
\left(
\sum_{r=1}^R \sum_{\bm{n}} \gamma_r(\bm{n}) \, \frac{1}{\bm{n}!} | \bm{n} + \bm{v}_r \rangle \langle \bm{n} |
\right)
| \bm{n}_{s} \rangle
=\sum_{r=1}^{R} \gamma_r(\bm{n}_s) \delta_{\bm{v}_r, \bm{n}_{s'}-\bm{n}_s}.
\end{eqnarray}
Hence, following tedious calculations based on the Gauss–Jordan elimination, the matrix element of \eref{eq_resolvent} is approximately expressed as:
\begin{eqnarray}
\fl
\left[ \rme^{\mathcal{L}^\dagger T} \right]_{s_\mathrm{f} s_\mathrm{i}} 
\simeq
\left( 1 - \frac{T}{M} \left[ \mathcal{L}^\dagger \right]_{s_\mathrm{i} s_\mathrm{i}} \right)
\left( \sum_{s_2} \cdots \sum_{s_{M-1}} \mathcal{R}_{s_\mathrm{i} s_{2}} \cdots \mathcal{R}_{s_{M-1} s_\mathrm{f}} \right)
\left( 1 - \frac{T}{M} \left[ \mathcal{L}^\dagger \right]_{s_\mathrm{f} s_\mathrm{f}} \right)^{-1} 
\label{eq_resolvent_element}
\end{eqnarray}
for large $M$, where 
\begin{eqnarray}
\left\{ \begin{array}{l}
\displaystyle \mathcal{R}_{ss}
= \frac{1}{1-\frac{T}{M} \left[ \mathcal{L}^\dagger \right]_{ss}},\\
\\
\displaystyle
\mathcal{R}_{s's} 
= \frac{T}{M} \left[ \mathcal{L}^\dagger \right]_{s' s}
\left( \frac{1}{1-\frac{T}{M} \left[ \mathcal{L}^\dagger \right]_{ss}} \right)^2 \quad \textrm{for } s' \neq s.
\end{array}
\right.
\end{eqnarray}
To derive \eref{eq_resolvent_element}, we assume that $(T/M)^2 \ll T/M $. A supplementary explanation for the derivation of \eref{eq_resolvent_element} is given in \ref{sec_appendix_inverse}. Repeated evaluation of an element of $\exp(\mathcal{L}^\dagger T)$ is conducted by increasing $M$, and the converged value provides the corresponding element of the Koopman matrix. Furthermore, when $M$ is not so large, the enumeration of the paths in \eref{eq_resolvent_element} is not difficult. Hence, a conventional dynamic programming method is available to evaluate \eref{eq_resolvent_element}.

Although the approximation for the resolvent in \eref{eq_resolvent_element} may be sufficient, it exhibits slow convergence behavior with an increase in $M$. Thus, a technique to avoid slow convergence is introduced here, that is, the extrapolation of sequences. There are various extrapolation methods; the following procedure is employed here:
\begin{enumerate}
\item Fix initial state $\bm{n}_{s_\mathrm{i}}$ and final state $\bm{n}_{s_\mathrm{f}}$. The indices $s_\mathrm{i}$ and $s_\mathrm{f}$ in \eref{eq_resolvent_element} determine the matrix element which we desire to evaluate.
\item The matrix elements \eref{eq_resolvent_element} are computed for $M = M_\mathrm{min}+1,M_\mathrm{min}+2,\dots,M_\mathrm{max}$. Let $\{\zeta_M^{(0)}\}$ be a sequence of the evaluated values. 
\item Consider the sequence $\{\zeta_M^{(0)}\}$ as a function of $1/M$.
\item Consider the finite differences for $\{\zeta_M^{(0)}\}$, and determine the tangent line at each point. Subsequently, the sequence of the intercept for each tangent line is obtained.
\item Repeat the previous procedure, and obtain the following sequences:
\begin{eqnarray}
\zeta_M^{(i+1)} = \zeta_M^{(i)} - \frac{\zeta_M^{(i)} - \zeta_{M-1}^{(i)}}{1/M - 1/(M-1)} \times \frac{1}{M}
\label{eq_diff_sequence}
\end{eqnarray}
for $i = 0, 1, \dots, N_\mathrm{diff}-1$; $N_\mathrm{diff}$ is the maximum degree of the finite differences.
\item Evaluate the weights
\begin{eqnarray}
w^{(i)} = \frac{1}{\left| \zeta_{M_\mathrm{max}}^{(i)} - \zeta_{M_\mathrm{max}-1}^{(i)} \right| + \varepsilon},
\end{eqnarray}
for $i=1,2,\dots,N_\mathrm{diff}$. In section~\ref{sec_numerical}, $\varepsilon = 10^{-5}$ is used. 
\item Consider the weighted sum as the final result:
\begin{eqnarray}
\left[ \rme^{\mathcal{L}^\dagger T} \right]_{s_\mathrm{f} s_\mathrm{i}} 
 \simeq \frac{w^{(i)}}{\sum_{i=1}^{N_\mathrm{diff}} w^{(i)}} \zeta_{M_\mathrm{max}}^{(i)}.
\label{eq_weighted_sum}
\end{eqnarray}
\end{enumerate}
Concrete examples of the sequences are provided in the next section, which gives a better understanding of this procedure.

There is a notice for the procedure. When we take the extrapolation, the case with $i=0$ is excluded here. Although it may be possible to include it, preliminary studies showed that its convergence for $i=0$ is worse than the others. Hence, we avoid it here.

\section{Numerical results}
\label{sec_numerical}

\subsection{Example: Noisy van der Pol process}
\label{subsec_van_der_Pol}

\begin{figure}
\begin{center}
\includegraphics[width=90mm]{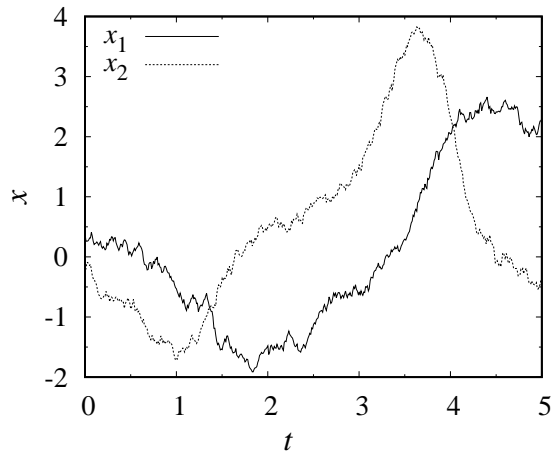}
\end{center}
\caption{A sample path for the noisy van der Pol system. The oscillating behavior is one of the features of the system.}
\label{fig_sample_path}
\end{figure}

To demonstrate the proposed method in section~\ref{sec_proposal}, the following example is considered:
\begin{eqnarray}
\left\{ \begin{array}{l}
\rmd x_1 = x_2 \rmd t + \nu_{11} \rmd W_1(t), \\
\rmd x_2 = \left( \epsilon x_2 \left( 1 - x_1^2 \right) - x_1\right) \rmd t + \nu_{22} \rmd W_2(t),
\end{array}\right.
\end{eqnarray}
which is the noisy version of van der Pol system \cite{van_der_Pol1926}. This stochastic process was previously used in the context of the Koopman operators; for example, see \cite{Crnjaric-Zic2020}. Here, we set $\epsilon = 1.0$, $\nu_{11} = 0.5$, and $\nu_{22} = 0.5$. Figure~\ref{fig_sample_path} shows sample paths.

The time-evolution operator for the dual process is expressed by
\begin{eqnarray}
\mathcal{L}^\dagger = x_2 \partial_{x_1} + \left( \epsilon x_2 \left( 1 - x_1^2 \right) - x_1\right) \partial_{x_2}
+ \frac{1}{2} \nu_{11}^2 \partial_{x_1}^2 + \frac{1}{2} \nu_{22}^2 \partial_{x_2}^2.
\end{eqnarray}
Furthermore, the shift of the origins is employed for $x_1$ and $x_2$ as follows:
\begin{eqnarray}
\fl
\mathcal{L}^\dagger = \left(x_2+x_2^\mathrm{c} \right) \partial_{x_1} 
+ \left( \epsilon \left(x_2+x_2^\mathrm{c} \right) 
\left( 1 - \left(x_1 + x_1^\mathrm{c} \right)^2 \right) - \left( x_1 + x_1^\mathrm{c} \right)\right) \partial_{x_2}
+ \frac{1}{2} \nu_{11}^2 \partial_{x_1}^2 + \frac{1}{2} \nu_{22}^2 \partial_{x_2}^2. \nonumber \\
\end{eqnarray}
To improve the numerical convergence, the shifting of the origins is sometimes beneficial. $\gamma_r(\bm{n})$ and $\bm{v}_r$ in the expression \eref{eq_simple_time_evolution_operator} are presented in Table 1.

\begin{table}
\caption{Coefficients and state-change vectors for the noisy van der Pol system.}
\begin{center}
\begin{tabular}{lll}
\hline\hline
Term & $\displaystyle \gamma_r(\bm{n})$ & $\displaystyle \bm{v}_r$\\
\hline
$(1/2) \nu_{11}^2 \partial_{x_1}^{2}$ 
& $ (1/2) \nu_{11}^2 n_1(n_1-1)$
& $[-2,0]$ \\
$x_{2}^{\mathrm{c}} \partial_{x_{1}}$ 
& $x_{2}^{\mathrm{c}} n_1$
& $[-1,0]$ \\
$x_{2} \partial_{x_{1}}$ 
& $n_1$
& $[-1,1]$ \\
$(1/2) \nu_{22}^{2} \partial_{x_{2}}^{2}$ 
& $(1/2) \nu_{22}^{2} n_2(n_2-1)$ 
& $[0,-2]$ \\
$- x_{1}^{\mathrm{c}} \partial_{x_{2}} - \epsilon \left(x_{1}^{\mathrm{c}}\right)^{2} x_{2}^{\mathrm{c}} \partial_{x_{2}} + \epsilon x_{2}^{\mathrm{c}} \partial_{x_{2}} \quad$ 
& $- x_{1}^{\mathrm{c}} n_2 - \epsilon \left(x_{1}^{\mathrm{c}}\right)^{2} x_{2}^{\mathrm{c}} n_2 + \epsilon x_{2}^{\mathrm{c}} n_2 \quad$
& $[0,-1]$ \\
$- \epsilon \left(x_{1}^{\mathrm{c}}\right)^{2} x_{2} \partial_{x_{2}} + \epsilon x_{2} \partial_{x_{2}} $
& $- \epsilon \left(x_{1}^{\mathrm{c}}\right)^{2} n_2 + \epsilon n_2$
& $[0,0]$\\
$- 2 \epsilon x_{1}^{\mathrm{c}} x_{2}^{\mathrm{c}} x_{1} \partial_{x_{2}} - x_{1} \partial_{x_{2}}$
& $- 2 \epsilon x_{1}^{\mathrm{c}} x_{2}^{\mathrm{c}} n_2 - n_2$
& $[1,-1]$\\
$- 2 \epsilon x_{1}^{\mathrm{c}} x_{1} x_{2} \partial_{x_{2}}$
& $- 2 \epsilon x_{1}^{\mathrm{c}} n_2$
& $[1,0]$\\
$ - \epsilon x_{2}^{\mathrm{c}} x_{1}^{2} \partial_{x_{2}} $
& $ - \epsilon x_{2}^{\mathrm{c}} n_2$
& $[2,-1]$\\
$- \epsilon x_{1}^{2} x_{2} \partial_{x_{2}} $
& $- \epsilon n_2$
& $[2,0]$\\
\hline\hline
\end{tabular}
\end{center}
\end{table}

\subsection{Results for the extrapolation}

In section~\ref{subsec_resolvent}, the extrapolation procedure was explained. Here, an example of the procedure is shown.

\begin{figure}
\begin{center}
\includegraphics[width=75mm]{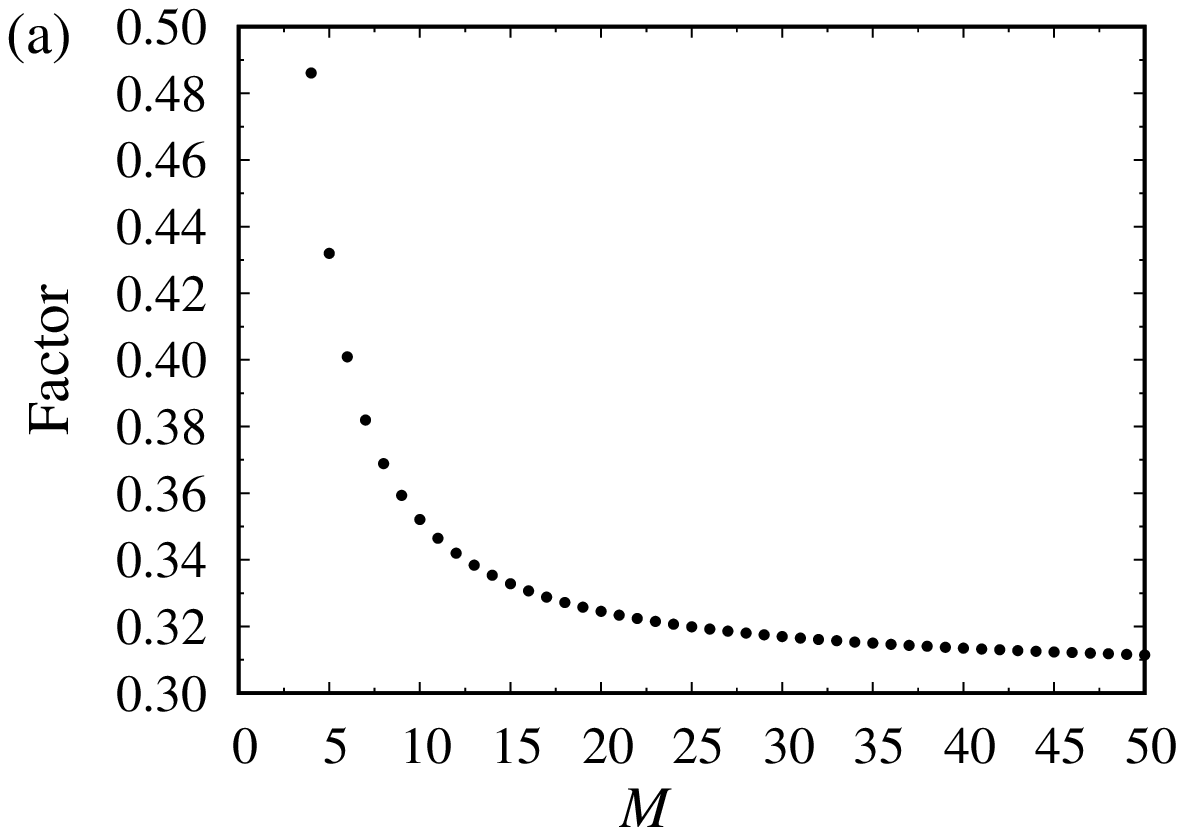}
\includegraphics[width=75mm]{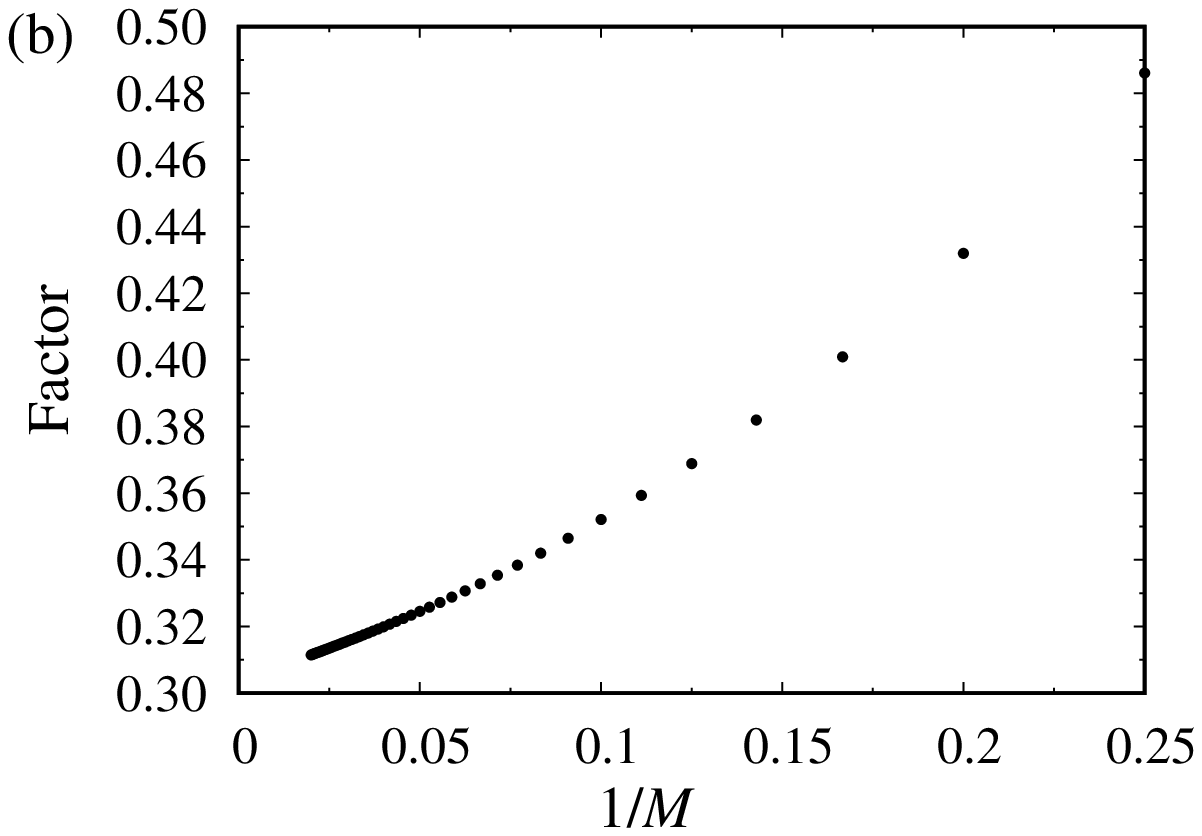}
\includegraphics[width=75mm]{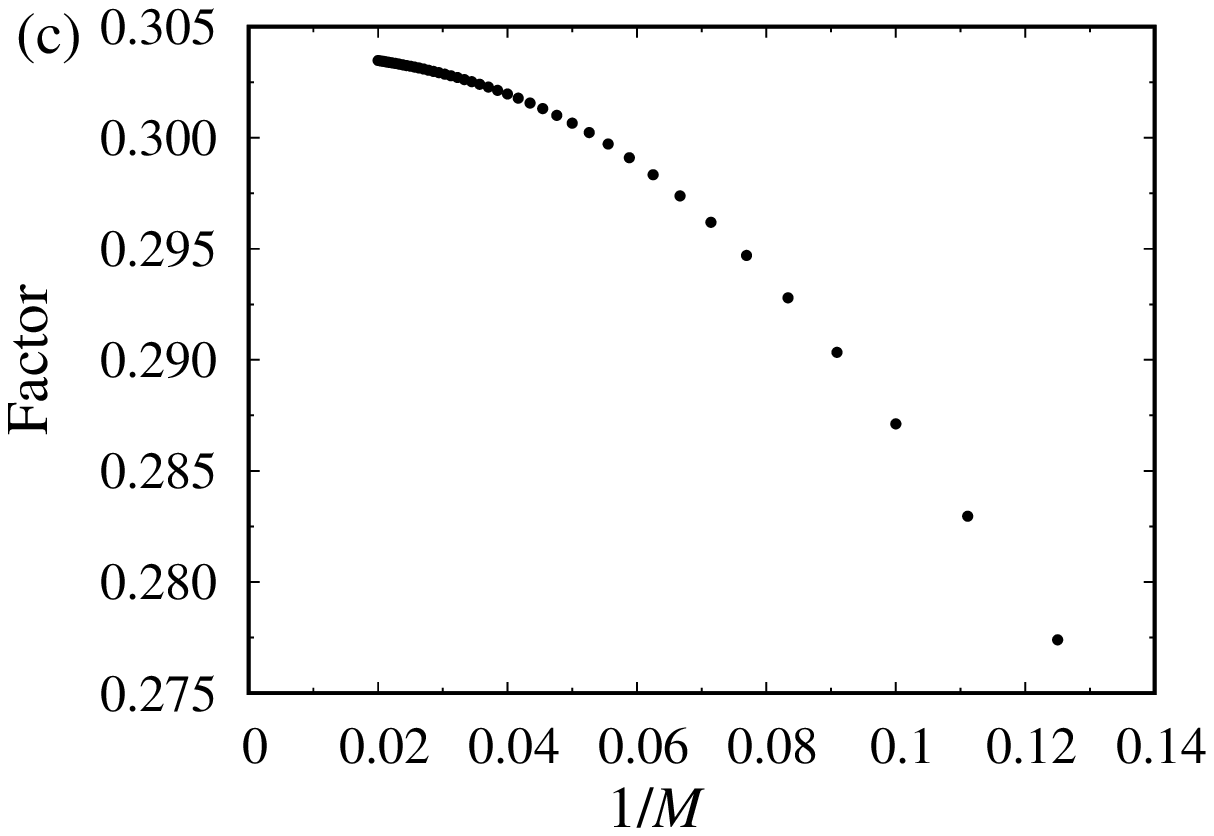}
\includegraphics[width=75mm]{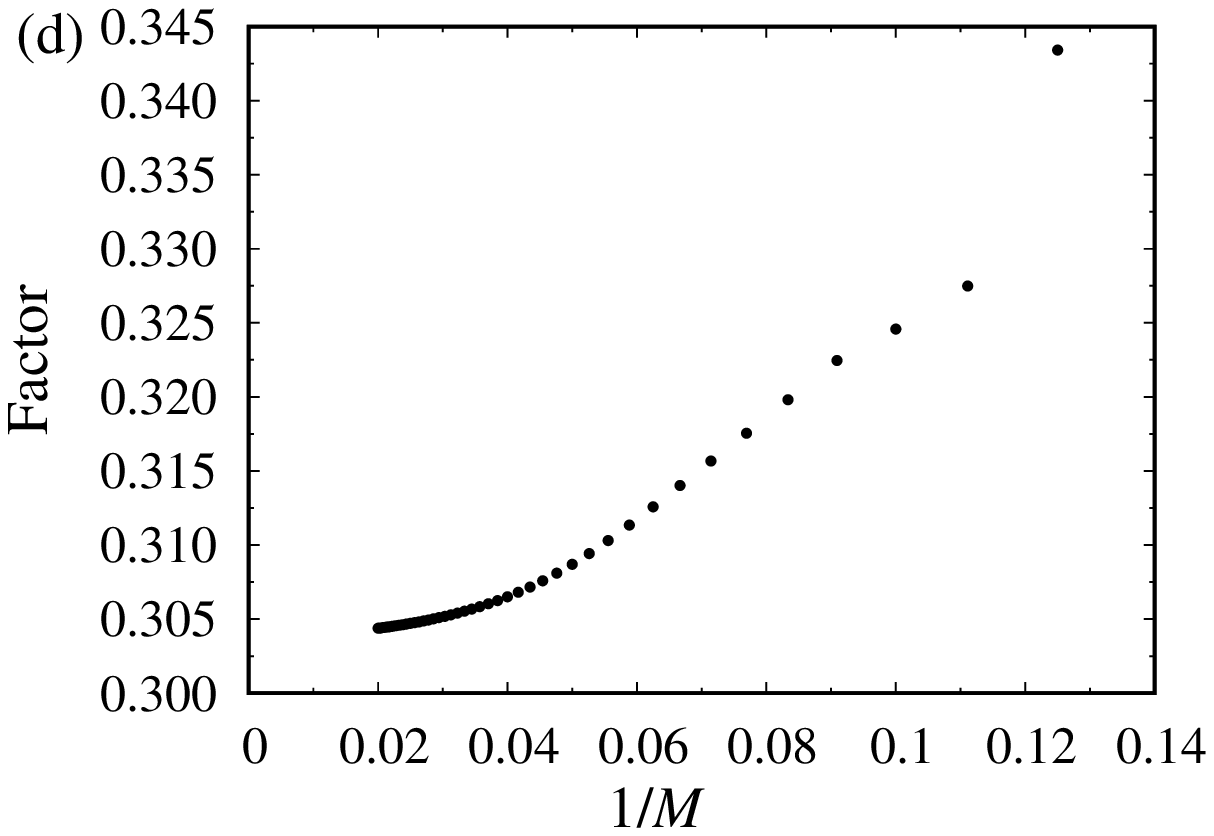}
\end{center}
\caption{Obtained sequences from the numerical evaluation of the resolvent. (a) The original sequence as the function of $M$. (b) The original sequence as the function of $1/M$, that is, $\{\zeta_M^{(0)}\}$. (c) The sequence $\{\zeta_M^{(1)}\}$ in \eref{eq_diff_sequence}. (d) The sequence $\{\zeta_M^{(2)}\}$ in \eref{eq_diff_sequence}.}
\label{fig_result_resolvent}
\end{figure}

Let $T = 1.0$, $x_1^\mathrm{c} = 0.0$ and $x_2^\mathrm{c} = 0.0$. Here, we calculate the Koopman matrix element for row $\langle \bm{n}_{s_\mathrm{f}} = (0,0)|$ and column $| \bm{n}_{s_\mathrm{i}} = (2,0) \rangle$ as a specific example; this element corresponds to evaluating $\mathbb{E}[x_1^2]$. The obtained sequence is plotted as a function of $M$, as shown in Figure~\ref{fig_result_resolvent}(a). The slow convergence of the sequence is evident. In contrast, when plotted as a function of $1/M$, the extrapolation at $1/M \to 0$ can be expected to yield good results, as shown in Figure~\ref{fig_result_resolvent}(b). Figures~\ref{fig_result_resolvent}(c) and (d) show the sequences $\{\zeta_M^{(1)}\}$ and $\{\zeta_M^{(2)}\}$, respectively. It should be noted that the scale of the vertical axis for each figure is smaller than that of Figure~\ref{fig_result_resolvent}(a). Based on these sequences, the Koopman matrix is evaluated without data.

Here, we consider $M_\mathrm{max} = 50$ and $N_\mathrm{diff} = 3$. The proposed method in section~\ref{sec_proposal} yields $3.0373 \times 10^{-1}$ as $s_\mathrm{f}s_\mathrm{s}$ element of the Koopman matrix. It is also possible to compute the element using Monte Carlo simulations because it corresponds to $\mathbb{E}[x_1^2]$. Hence, the value of $\mathbb{E}[x_1^2]$ was estimated by the Monte Carlo simulations here by using $1000$ samples generated with the Euler-Maruyama approximation with the small-time interval $10^{-3}$ \cite{Kloeden_book}. In addition, the same Monte Carlo simulations were repeated $100$ times to obtain the average of the expectation value and standard deviation. Finally, $3.0451 \times 10^{-1} \pm 1.40 \times 10^{-2}$ was obtained as the numerical estimation via the Monte Carlo simulations. Thus, the proposed method is confirmed to work well.

There is a comment regarding the computational time. The above Monte Carlo procedure requires approximately $7$ seconds in total. In contrast, the proposed method requires approximately $0.15$ seconds. Although the computational time of the proposed method depends on the parameters such as $M_\mathrm{max}$, it is clear that the proposed method is beneficial from the perspective of computational time.

\subsection{Comparison with the EDMD method}

\begin{figure}
\begin{center}
\includegraphics[width=110mm]{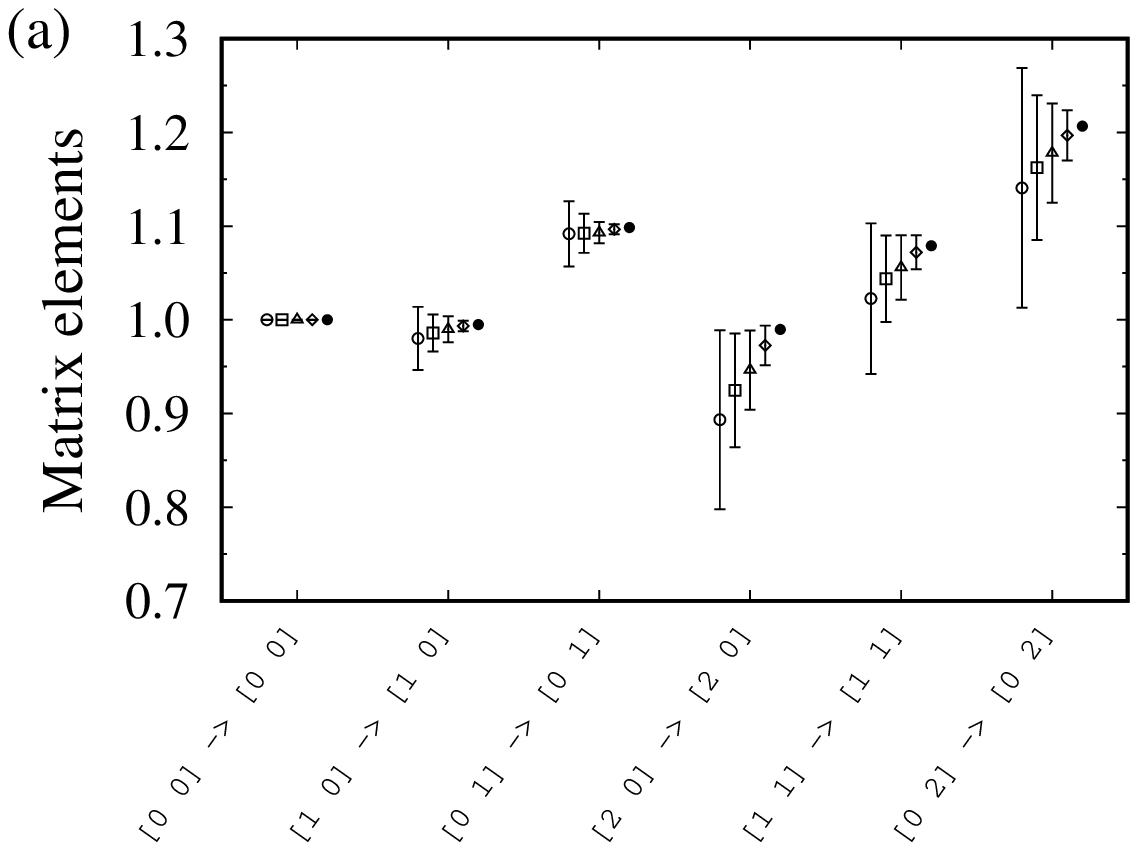}
\includegraphics[width=110mm]{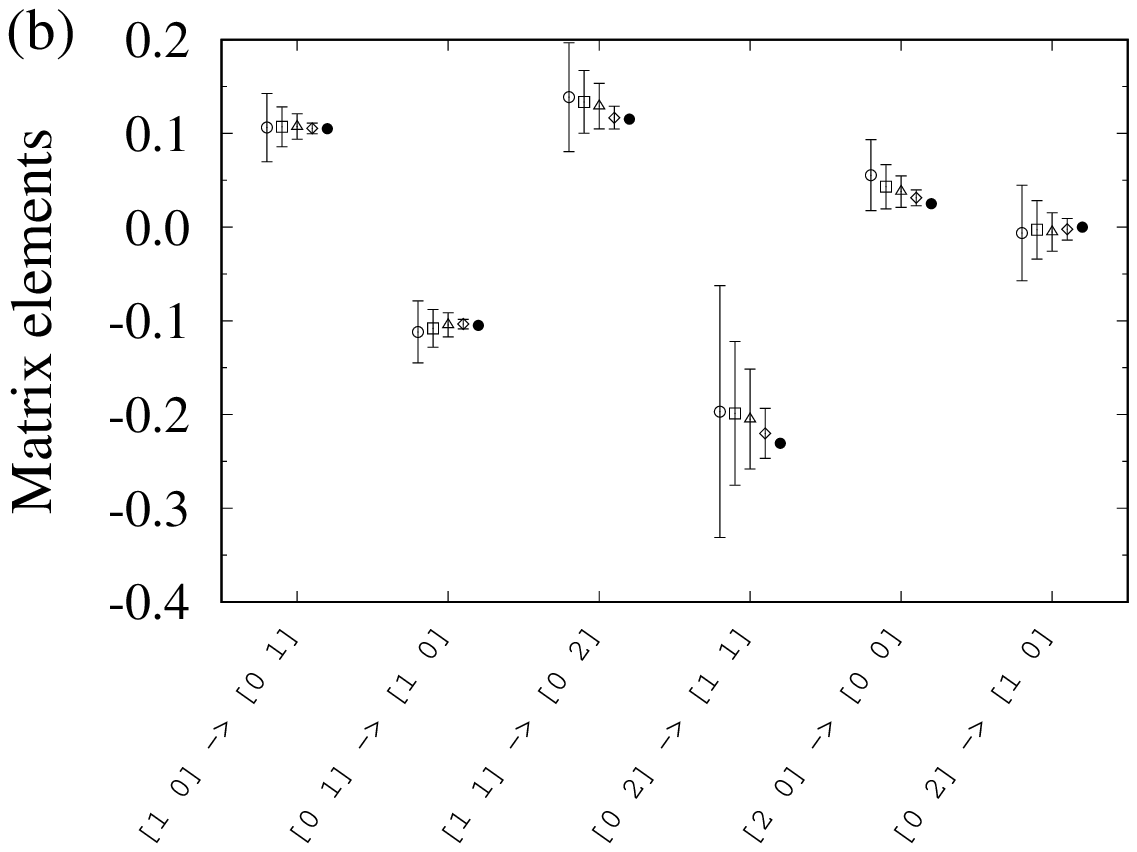}
\end{center}
\caption{Estimations of certain Koopman matrix elements. The circles, boxes, triangles, and diamonds correspond to the estimates for data size $N_{10} = 500$, $N_{20} = 1000$, $N_{40} = 2000$, and $N_{160} = 8000$, respectively. The estimates obtained from the proposed method are depicted with the filled circles. Results for several diagonal elements are shown in (a), and those for non-diagonal ones are shown in (b).}
\label{fig_result_Koopman_matrix}
\end{figure}

Finally, a comparison between the EDMD method and the proposed procedure is conducted. Herein, the following parameters are used: $x_1^\mathrm{c} = 0.0, x_2^\mathrm{c} = 0.0, T = 0.1, M_\mathrm{max} = 50$ and $N_\mathrm{diff} = 3$. For the EDMD, we employ the Euler-Maruyama method with the small time-interval $10^{-3}$, and sample paths with $t_\mathrm{end} = 5.0$ are generated. Because we set $T = 0.1$, it is possible to obtain $50$ snapshot pairs from a sample path. The dictionary for the EDMD is chosen as the monomial-type, with the largest degree set as $6$th. The EDMD method is applied $100$ times, and thereafter the average and standard deviation were evaluated.

Figure~\ref{fig_result_Koopman_matrix} shows the result of the comparisons. The error bars are standard derivations. The horizontal axis shows the element of the Koopman matrix; for example, the description \verb|[1 0] -> [0 1]| implies $\bm{n}_{s_\mathrm{i}} = (1,0)$ and $\bm{n}_{s_\mathrm{f}} = (0,1)$. When $10$ sample paths are used for the EDMD, the number of snapshot pairs is $500$, as explained above. Let $N_{10}$ be the size of data for $10$ sample paths. Then, $N_{10} = 500$. Cases with different data size, $N_{20} = 1000$, $N_{40} = 2000$, and $N_{160} = 8000$ are numerically simulated. The circles, boxes, triangles, and diamonds in Figures~\ref{fig_result_Koopman_matrix}(a) and (b) correspond to $N_{10}, N_{20}, N_{40}$, and $N_{160}$ cases, respectively. The filled circles in Figures~\ref{fig_result_Koopman_matrix}(a) and (b) indicate the results of the proposed procedure. Figure~\ref{fig_result_Koopman_matrix}(a) shows the diagonal cases; Figure~\ref{fig_result_Koopman_matrix}(b) shows the non-diagonal cases. Note that the speeds of convergence in the EDMD are different for the elements. Here, the \verb|[0 2] -> [1 0]| case is exactly zero, and the proposed procedure provides the adequate result. In contrast, the EDMD exhibits large errors for small data size.

Although only certain typical cases are illustrated in Figure~\ref{fig_result_Koopman_matrix}, the numerical experiments show that the proposed procedure works well. 

Numerical results for another example is shown in \ref{sec_appendix_OU}; the proposed method gives reasonable results for the Ornstein-Uhlenbeck process.

\section{Discussion and conclusion}

The Koopman matrix is linear and tractable; its role would become increasingly important. There are many recent studies for the data-driven methods for the Koopman operators. The present work focused on the numerical evaluation of the Koopman matrix from the system equations; the proposed algorithm can evaluate directly only some of the elements of the Koopman matrix; there is no need to calculate the whole matrix at once. The proposed method mitigates the computational difficulties, in which the usage of resolvent and extrapolation also contribute to the reduction of the computational time. 

There are two complementary approaches to the Koopman matrix. One is the conventional data-driven method, and another is the proposed approach based on system equations. The combination of these two approaches would give more efficient numerical frameworks for evaluating the Koopman matrix. For example, firstly, we perform a rough estimation for the Koopman matrix from prior knowledge of system equations, and then, the rough estimation is modified with small-size data. The proposed method in this study will contribute to future research to construct such an integrated approach.

There are some remarks for the present work. This work is the first attempt to calculate the Koopman matrix with the numerical algorithm based on the resolvent. One of the remarkable features is the direct evaluation of an element of the Koopman matrix, as stated above. This feature would be beneficial for high-dimensional cases to avoid the curse of dimensionality. Of course, further mathematical and numerical works are necessary for the future. For example, the evaluation of the resolvent should be performed on the resolvent set, and hence we should check that the proposed numerical procedure does not violate this condition. In preliminary numerical studies, the proposed algorithm works well for several examples if $T$ is not so large; for example, see \cite{Ohkubo2021} for the research in the context of the duality. However, it will be necessary to check the condition with mathematically rigorous discussions in future works. Another question is the stability of the algorithm. The adjoint operator $\mathcal{L}^\dagger$ is usually unbounded, and hence a conventional Taylor-type expansion for the exponential function does not work well. As stated in \cite{Kato_book}, the usage of resolvents can avoid this problem. However, it is preferable to discuss the stability of the proposed method from the viewpoint of mathematics. In functional theory, it is known that the norm of a strongly continuous semigroup is bounded for a fixed time \cite{Engel_book}; i.e., the exponential of $\mathcal{L}^\dagger T$ is bounded when we fix $T$. The viewpoint of operator semigroups could be useful for the discussion. In addition, the choice of the basis function will affect the performance of the numerical algorithm. In the context of duality, there are some discussions about the choice of orthogonal polynomials \cite{Franceschini2019}. It should be important to discuss the properties of the chosen dictionaries for the Koopman operator cases.

As stated above, there are several remaining works from the viewpoints of mathematics and algorithms. These are beyond the scope of this work. However, the demonstration given in this work opens a new way for the numerical evaluation of the Koopman matrix. The present work will stimulate other researchers in mathematics and algorithms.

\ack
This work was partially supported by JST, PRESTO Grant Number JPMJPR18M4, Japan, and by JSPS KAKENHI Grant Number JP21K12045.
\vspace{3mm}

\appendix

\section{Supplementary explanation for \eref{eq_path_integral_like}}
\label{sec_appendix_path_integral}

A simple Taylor expansion of the exponential function gives \eref{eq_path_integral_like}. That is,
\begin{eqnarray}
e^{\mathcal{L}^\dagger T} |\bm{n}_\mathrm{ini} \rangle 
&= \sum_{M=0}^{\infty} \frac{T^M}{M!} \left( \mathcal{L}^\dagger \right)^M |\bm{n}_\mathrm{ini} \rangle \nonumber \\
&= \sum_{M=0}^{\infty} \frac{T^M}{M!}
\left( \sum_{r=1}^R \sum_{\bm{n}} \gamma_r(\bm{n})   \frac{1}{\bm{n}!}  | \bm{n} + \bm{v}_r \rangle \langle \bm{n} | \right)^M |\bm{n}_\mathrm{ini} \rangle.
\end{eqnarray}
Note that 
\begin{eqnarray}
 \sum_{r=1}^R \sum_{\bm{n}} \gamma_{r}(\bm{n})   \frac{1}{\bm{n}!}  | \bm{n} + \bm{v}_{r} \rangle \langle \bm{n} |\bm{n}_\mathrm{ini} \rangle 
&= 
\sum_{r=1}^R \sum_{\bm{n}} \gamma_r(\bm{n}) \frac{1}{\bm{n}!}| \bm{n} + \bm{v}_r \rangle \bm{n}_{\mathrm{ini}}! \delta_{\bm{n},\bm{n}_{\mathrm{ini}}} \nonumber \\
&=
 \sum_{r=1}^R \gamma_r(\bm{n}_{\mathrm{ini}})  | \bm{n}_{\mathrm{ini}} + \bm{v}_r \rangle,
\label{eq_appendix_path_integral_first}
\end{eqnarray}
where \eref{eq_bra_ket_orthogonal} is used. Then, the same procedure is repeated $M$ times. Finally, the derived terms are rearranged by focusing the number of events $M$, and \eref{eq_path_integral_like} is derived.

\section{Supplementary explanation for \eref{eq_resolvent_element}}
\label{sec_appendix_inverse}

It would be beneficial to consider a simple $2 \times 2$ matrix and its inverse matrix. The Gauss-Jordan elimination \cite{Strang_book} gives the following transformation using elementary row operations:
\begin{eqnarray}
&\left(
\left.
\begin{array}{cc}
a_{11} & a_{12} \\
a_{21} & a_{22} 
\end{array}
\right|
\begin{array}{cc}
1 & 0 \\
0 & 1 
\end{array}
\right)
\to
\left(
\left.
\begin{array}{cc}
1 & a_{12}/a_{11} \\
a_{21} & a_{22}
\end{array}
\right|
\begin{array}{cc}
1 / a_{11} & 0 \\
0 & 1
\end{array}
\right) \nonumber \\
&\to
\left(
\left.
\begin{array}{cc}
1 & a_{12}/a_{11} \\
0 & a_{22} - (a_{21}a_{12})/a_{11}
\end{array}
\right|
\begin{array}{cc}
1 / a_{11} & 0 \\
- a_{21}/a_{11} & 1
\end{array}
\right). \nonumber \\
\end{eqnarray}
We here assume that $a_{22} \gg (a_{21}a_{12}/a_{11})$, and an approximation of $a_{22} - (a_{21}a_{12}/a_{11}) \simeq a_{22}$ is employed. Hence,
\begin{eqnarray}
\left(
\left.
\begin{array}{cc}
1 & a_{12}/a_{11} \\
0 & a_{22} 
\end{array}
\right|
\begin{array}{cc}
1 / a_{11} & 0 \\
- a_{21}/a_{11}  & 1
\end{array}
\right)
\to
\left(
\left.
\begin{array}{cc}
1 & a_{12}/a_{11} \\
0 & 1
\end{array}
\right|
\begin{array}{cc}
1 / a_{11} & 0 \\
- a_{21}/(a_{11}a_{22})  & 1/a_{22}
\end{array}
\right) \nonumber \\
\to
\left(
\left.
\begin{array}{cc}
1 & 0 \\
0 & 1
\end{array}
\right|
\begin{array}{cc}
1 / a_{11} + a_{21}/(a_{11}a_{22})\times a_{12}/a_{11} &  - a_{12}/(a_{11} a_{22}) \\
- a_{21}/(a_{11}a_{22})  & 1/a_{22}
\end{array}
\right). \nonumber \\
\end{eqnarray}
After taking an approximation $1 / a_{11} + a_{21}/(a_{11}a_{22})\times a_{12}/a_{11} \simeq 1/a_{11}$, we finally obtain
\begin{eqnarray}
A^{-1} = \left(
\begin{array}{cc}
1 / a_{11}  &  - a_{12}/a_{11} a_{22} \\
- a_{21}/a_{11}a_{22}  & 1/a_{22}
\end{array}
\right).
\end{eqnarray}
Applying the similar discussion for $1 - (T/M)\mathcal{L}^\dagger$ and assuming $T/M \ll 1$, the approximations could be verified and we obtain
\begin{eqnarray}
\fl
\left\{ \begin{array}{l}
\displaystyle \left[ 1 - \frac{T}{M} \mathcal{L}^\dagger \right]_{ss}^{-1}
= \frac{1}{1-\frac{T}{M} \left[ \mathcal{L}^\dagger \right]_{ss}}\\
\\
\displaystyle
\left[ 1 - \frac{T}{M} \mathcal{L}^\dagger \right]_{s's}^{-1}
= \frac{T}{M} \left[ \mathcal{L}^\dagger \right]_{s' s}
\left( \frac{1}{1-\frac{T}{M} \left[ \mathcal{L}^\dagger \right]_{s's'}} \right)
\left( \frac{1}{1-\frac{T}{M} \left[ \mathcal{L}^\dagger \right]_{ss}} \right)
 \quad \textrm{for } s' \neq s.
\end{array}
\right.
\label{eq_appendix_inverse}
\end{eqnarray}

Since the subscript $s' s$ means the transition $s \to s'$, it is possible to rearrange the denominators. Finally, \eref{eq_resolvent_element} is derived after multiplying the complementary factors related to the initial and final states. From the viewpoint of computational algorithms, the expression in \eref{eq_resolvent_element} reduces the number of memory accesses in the numerical calculation for elements of the inverse matrix compared with the original expression \eref{eq_appendix_inverse}.

\section{Numerical results for Ornstein-Uhlenbeck process}
\label{sec_appendix_OU}

Although the noisy van der Pol process is used in section~\ref{sec_numerical}, it would be helpful to show simpler example cases. The following Ornstein-Uhlenbeck process is a good candidate:
\begin{eqnarray}
dx = - \gamma x dt + \sigma dW(t),
\end{eqnarray}
where $\gamma \in \mathbb{R}$ and $\sigma \in \mathbb{R}_{+}$. It is possible to obtain analytical solutions for the Koopman eigenfunctions. For example, \cite{Crnjaric-Zic2020} gives a numerical study for the Koopman eigenfunction for the Ornstein-Uhlenbeck process with an algorithm based on the dynamic mode decomposition. Here, we aim to evaluate elements of the Koopman matrix. Hence, the same procedure in section~\ref{sec_numerical} is employed. Comparisons between the numerical results of the EDMD and those of the proposed method are given.

A remarkable fact for the Ornstein-Uhlenbeck process is that the calculation of the resolvent is quite easy because of the non-increasing feature of the process. That is, the adjoint operator for the original process is given as
\begin{eqnarray}
\mathcal{L}^\dagger = -\gamma x \partial_x + \frac{\sigma^2}{2} \partial_x^2.
\label{eq_appendix_OU_operator}
\end{eqnarray}
The first term in the r.h.s. in \eref{eq_appendix_OU_operator} does not change the state $n$, and the second term acts as $n \to n-2$. Hence, there is no transition to increase the state. In the dynamic programming method for solving \eref{eq_resolvent_element}, the number of states generally increases as $M$ increases. However, it is enough to consider only a small number of states for the Ornstein-Uhlenbeck process, which is preferable for numerical calculations.

\begin{figure}
\begin{center}
\includegraphics[width=110mm]{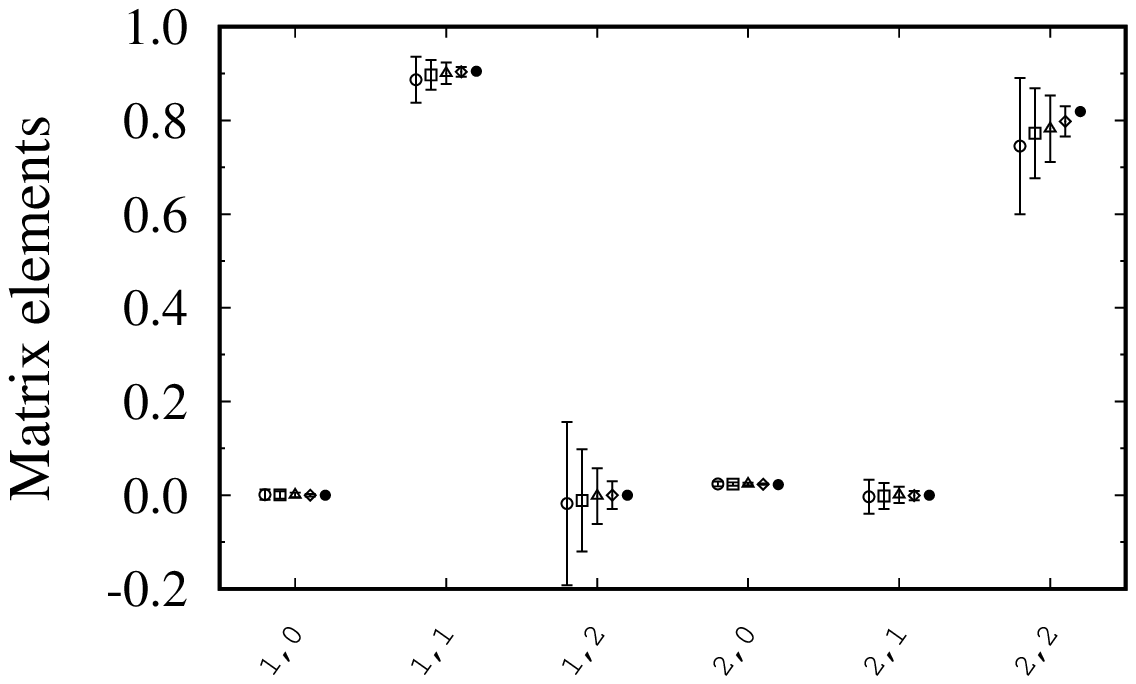}
\vspace{-15mm}
\end{center}
\caption{
Estimations of certain Koopman matrix elements for the Ornstein-Uhlenbeck process. The circles, boxes, triangles, and diamonds correspond to the estimates for data size $N_{10} = 500$, $N_{20} = 1000$, $N_{40} = 2000$, and $N_{160} = 8000$, respectively. The estimates obtained from the proposed method are depicted with the filled circles.}
\label{fig_result_Koopman_matrix_OU}
\end{figure}

Figure~\ref{fig_result_Koopman_matrix_OU} shows the numerical results for $\gamma = 1.0$ and $\sigma = 0.5$. The other settings are the same as in section~\ref{sec_numerical}. Because of the non-increasing feature stated above, the matrix element corresponding to $1 \to 2$ is exactly zero. Although the EDMD has large estimation errors for the $1\to 2$ case, the proposed method adequately recovers this exact result.

\end{document}